\documentclass[a4paper]{scrartcl}
\RequirePackage{fix-cm}
\usepackage{graphicx}
\usepackage{mathptmx}      

\usepackage[utf8]{inputenc}
\usepackage{amssymb}
\usepackage{amsmath}
\usepackage{amsmath,amsthm}

\usepackage{algorithm}
\usepackage{algpseudocode}
\usepackage[colorlinks,citecolor=blue]{hyperref} 

\newtheorem{theorem}{Theorem}
\newcommand{\LL}{{\cal L}}
\newcommand{\RR}{{\mathbb R}}
\newcommand{\vv}[1]{#1} 
\newcommand{\bhat}{\hat{b}}
\newcommand{\basis}[2]{\vect{#1}{0},\vect{#1}{1},\ldots,\vect{#1}{#2-1}}
\newcommand{\basisfull}[2]{\vect{#1}{0},\vect{#1}{1},\ldots,\vect{#1}{#2}}
\newcommand{\awproj}[1]{\pi_{#1}}              
\newcommand{\ip}[2]{\langle #1,#2\rangle}  
\newcommand{\spann}[2]{\langle\basis{#1}{#2}\rangle}

\newcommand{\ZZ}{\mathbb{Z}}
\newcommand{\vect}[2]{{#1}^{(#2)}}
\newcommand{\rmax}{r_{\max}}
\DeclareMathOperator{\lcm}{lcm}

\begin{document}
\title{Solving the Market Split Problem with Lattice Enumeration}
\author{Alfred Wassermann\\
    Department of Mathematics, University of Bayreuth\\
    95447 Bayreuth, Germany \\
    \texttt{alfred.wassermann@uni-bayreuth.de}           
}

\maketitle

\begin{abstract}
    The market split problem was proposed by Cornu\'ejols and Dawande in 1998
    as benchmark problem for algorithms solving linear systems
    with binary variables.
    The recent (2025) Quantum Optimization Benchmark Library (QOBLIB)
    contains a set of feasible instances of the market split problem.
    In QOBLIB an instance of the market split problem is considered as solved
    as soon as at least one feasible solution has been found.

    The market split problem seems to be difficult to solve 
    with the conventional branch-and-cut approach of
    integer linear programming software which reportedly can handle 
    QOBLIB instances up to $m=7$. In contrast, a new GPU implementation of the
    Schroeppel--Shamir algorithm solves instances up to $m=11$.
    
    In this note we report about experiments
    with an algorithm that reduces the market split problem to a lattice problem.
    With the author's most recent implementation -- named \texttt{solvediophant} --
    instances of the QOBLIB market split benchmark problems
    can be solved up to $m=14$ on a standard computer.

\end{abstract}


\section{Introduction}\label{intro}
In~\cite{CornuejolsDawande:98,CornuejolsDawande:99}, Cornu\'ejols and Dawande proposed the \emph{market split problem} (MSP) as benchmark problem
for integer linear programming software.
Its formulation as optimization problem is:
\begin{eqnarray*}
    &\min & \sum_{i=1}^m |s_i| \\
    &\mbox{s.t.}&  \sum_{j=1}^n a_{ij}x_j +s_i = d_i, \qquad i=1,\ldots ,m\\
    && x_j \in\{0,1\}, \qquad j=1,\ldots,n \\
    && s_i\in\ZZ, \qquad i=1,\ldots ,m,\\
\end{eqnarray*}
The numbers can be interpreted like this:
$n$ is the number of retailers, $m$ is the number of products,
$a_{ij}$ is the demand of retailer $j$ for product $i$, and the right
hand side vector $d_i$ is determined from the desired market split
among the divisions $D_1$ and $D_2$ of a company, see~\cite{CornuejolsDawande:98,Williams:78}.

The feasibility version (fMSP) of the problem is:
\[
    \exists  \; x \in\{0,1\}^n
    \quad\mbox{s.t.}\;  \sum_{j=1}^n a_{ij}x_j = d_i,\quad i=1,\ldots ,m \; ?\\
\]
The problem is NP-complete, since for $m=1$ it is the subset sum problem
which is known to be NP-complete~\cite{GareyJohnson:79}.
Since these problems seem to be difficult to solve with
the conventional branch-and-cut approach even for relatively
small values of $m$ and $n$, the authors
of~\cite{CornuejolsDawande:98} presented them as a challenge
for the research community.

In~\cite{CornuejolsDawande:98}, the following
special instances are proposed as test problems, depending on a single parameter $m$:
\begin{quote}
    Fix $n =10\cdot(m-1)$.
    The numbers $a_{ij}\in \ZZ$ are generated
    uniformly and independently at random with $0 \leq a_{ij} < 100$,
    $1\leq i\leq m$, $1\leq j\leq n$,
    and the right hand side is set to $d_i = \lfloor\frac{1}{2}\sum_{j=1}^n a_{ij}\rfloor$.
\end{quote}
In~\cite[2002]{wassermann:02}, the present author
solved instances of (fMSP) up to $m=10$ with his software \texttt{solvediophant}.
See~\cite{KempkeKochZIB:25} for a list of publications about
the market split problem since 2002.

The recent
\emph{Quantum Optimization Benchmark Library (QOBLIB)}~\cite{koch2025quantumoptimizationbenchmarklibrary}
proposes a slight variation of
the original market split problem (fMSP) as one of its classes of benchmark problems,
i.e.~it contains problems of the form
\[
    n =10\cdot(m-1), \quad
    a_{ij}\in \ZZ, \;0 \leq a_{ij} < D,
\]
depending on two parameters $m$ and $D$.
For every combination of $m\in\{3, \ldots, 15\}$ and $D\in\{50, 100, 200\}$,
QOBLIB contains four instances of (fMSP), and it has been made sure that every instance
has at least one solution.
In the original formulation of the market split problem, an instance is \emph{solved} if the
feasibility question is answered. That means usually, the answer is ``yes'' as soon as
one solution for an instance has been found.
However, since in~\cite{koch2025quantumoptimizationbenchmarklibrary} all instances are feasible by construction and
the task for each instance is formulated as ``to find $x\in\{0,1\}^n$ such that $Ax=b$'',
here an instance of the market split problem is regarded as \emph{solved} as soon as one $x\in\{0,1\}^n$ has been found such that $Ax=b$.

In~\cite{KempkeKochZIB:25}, see also~\cite{KempkeKoch:25}, the Schroeppel--Shamir algorithm~\cite{SchroeppelShamir:1981}
was implemented on a GPU and used to solve QOBLIB market split instances successfully up to $m=11$.

\medskip
Here we report about experiments with the author's version 1.0
of \texttt{solvediophant} applied to the QOBLIB market split instances,
see Table~\ref{table:timings}. The program could solve all instances up to $m=14$ and
exhaustively find all binary solution vectors up to $m=9$.

\section{Algorithm}

The algorithm that was used to solve market split problems from
QOBLIB allows to solve the following, slightly more general problem:

Let $A\in \ZZ^{m\times n}$, $d\in \ZZ^m$, and
$l, r \in \ZZ^n$.
Determine all vectors $x \in \ZZ^n$ such that
\begin{equation}
    A\cdot x = d \;\mbox{ and }\; l \leq x \leq r\,, \label{diophant}
\end{equation}
where $l\leq r$ for vectors $l, r\in \ZZ^n$ is defined as
$l_i\leq r_i$ for all $0\leq i < n$.

It is worth to emphasize that for this algorithm 
the matrix $A$ and the right hand side vector $d$ might have negative entries, too.
Moreover, with the substitution $x:= x - l$, $d:= d- A\cdot l$
and $r:= r - l$, it suffices to consider $l = 0$ as a lower bound on the variables.

Here, we give only an informal description of the algorithm, for details
we refer to~\cite{Wassermann:98,wassermann:02,10.1007/978-3-030-79987-8_2}.
In order to simplify the description of the algorithm, we also assume that
$m \leq n$ and the rows of the matrix $A$ are linearly independent. In general, the
algorithm works without changes for matrices with an arbitrary number of rows and $m$ being
the rank of $A$ (over the rationals).

Problem~(\ref{diophant}) is reduced to a lattice problem by considering the lattice
spanned by the columns of the $(m+n+1)\times (n+1)$ matrix
\begin{equation}
    \left(
    \begin{array}{c|cccc}
            -N\cdot d & \multicolumn{4}{c}{N\cdot A}                            \\\hline
            -\rmax    & 2c_1                         & 0      & \cdots & 0      \\
            -\rmax    & 0                            & 2c_2   & \cdots & 0      \\
            \vdots    & \vdots                       &        & \ddots & \vdots \\
            -\rmax    & 0                            & \cdots & \cdots & 2c_n   \\\hline
            \rmax     & 0                            & \cdots & \cdots & 0      \\
        \end{array}
    \right),															\label{lattice}
\end{equation}
where $N\in\ZZ_{>0}$ is a large constant and
$$
    \rmax = \lcm\{r_1,\ldots,r_n\}
    \quad\text{and}\quad
    c_i = \frac{\rmax}{r_i}\;, \quad 1\leq i\leq n\,.
$$
In case of a problem with restriction $x\in\{0,1\}^n$, $\rmax = c_1=\cdots=c_n = 1$.

Theoretical lower bounds for $N$ are given in \cite[Thm 4]{AardalHurkensLenstra:98}.
However, in most cases lattice reduction performs better than predicted and allows $N$ to
be much smaller than these lower bounds.
In practice, for (fMSP) instances of QOBLIB, setting in \texttt{solvediophant} $N=2^{10}$ suffices.

In a first step of the algorithm, the lattice basis is reduced with
the LLL algorithm or blockwise Korkine--Zolotarev reduction,
see~\cite{10.1007/978-3-030-79987-8_2} and the references therein.
If $N$ is sufficiently large,
the reduced basis consists of $n-m+1$ vectors with only zeroes in the first $m$ rows and
$m$ vectors which contain at least one nonzero entry in the first $m$ rows.
The latter vectors can be removed from the basis. From the remaining $n-m+1$ vectors
we can delete the first $m$ rows which contain only zeroes.
This gives a basis
$\vect{b}{0}$, $\vect{b}{1}$, $\ldots$, $\vect{b}{n-m}\in\ZZ^{n+1}$ of the kernel of the extended
system $Ax - d x_{n+1} = 0$.

The second step of the algorithm exhaustively enumerates all integer linear combinations
of the basis vectors $\vect{b}{0}$, $\vect{b}{1}$, $\ldots$, $\vect{b}{n-m}\in\ZZ^{n+1}$
which correspond to solutions of \eqref{diophant} as stated in the following theorem.

\begin{theorem}[\cite{wassermann:02}]
    Let
    \begin{equation}
        w = u_0\cdot \vect{b}{0}+u_1\cdot \vect{b}{1}+
        \ldots+u_{n-m}\cdot \vect{b}{n-m}        \label{probleminfinite}
    \end{equation}
    be an integer linear combination of the basis vectors
    with $w_0=\rmax$.
    $w$ is a solution of the system~\eqref{diophant} if and only if
    \[
        w\in\ZZ^{n+1} \mbox{ where } -\rmax\leq w_i \leq \rmax, \; 1\leq i\leq n\;.
    \]
\end{theorem}

The overall approach is related to~\cite{AardalBixbyHurkens:99,AardalHurkensLenstra:98},
which use lattice basis reduction too, but with a different lattice
and different exhaustive enumeration. A further distinction of the algorithm described here is
that it uses H\"older's inequality in the enumeration part of the algorithm. This is an
idea that has been proposed by~\cite{Ritter:97} for $\{0,1\}$ problems.
Here we solve the general problem with arbitrary bounds on the variables.



A priori, a lattice $\LL=\{\sum_{i=0}^{n-m} u_i\vect{b}{i}\mid u_i\in\ZZ\}$
of rank $n-m+1$ contains infinitely many elements.
However, it will turn out that the there are bounds on the integers
$|u_i|$, $0\leq i< n-m$ which lead to solutions of \eqref{diophant}, which depend solely on the lattice basis
$\vect{b}{0}$, $\vect{b}{1}$, $\ldots$, $\vect{b}{n-m}$.

These bounds reduce the problem of finding solution vectors to a finite set of lattice vectors.
Each solution vector $\vv v$ has the upper bounds
\[
    \|\vv v\|_2^2 \leq (n+1)\cdot\rmax^2 \quad \mbox{ and } \quad \|\vv v\|_\infty \leq \rmax\,.
\]
The exhaustive enumeration is arranged as backtracking algorithm.
Starting from $u_{n-m}\in\ZZ$, successively
all possible $u_t\in\ZZ$ for $t=n-m,n-m-1,\ldots,1,0$ are tested.
The enumeration can be pruned at stage $t$ if certain conditions are violated. These pruning tests
have quite a long history and are based on the work
of~\cite{CoveyouMacPherson:67,Dieter:75,FinckePohst:85,KaibRitter:95,Kannan:87,Knuth:69,Ritter:97}.

We recall that the \emph{Gram--Schmidt orthogonalization} of a basis $\vect{b}{0}$, $\vect{b}{1}$, $\ldots$, $\vect{b}{n-m}$
is the orthogonal family $\vect{\bhat}{i},\, 0\leq i \leq n-m$ defined by
\begin{equation}
    \vect{\bhat}{i} = \vect{b}{i} - \sum_{j=0}^{i-1}\mu_{ij}\cdot \vect{\bhat}{j}\, ,
    \;\mbox{where}\;
    \mu_{ij} =
    \frac{\ip{\vect{b}{i}}{\vect{\bhat}{j}}}{\ip{\vect{\bhat}{j}}{\vect{\bhat}{j}}}
    \label{gramschmidt}
    \;\mbox{and}\;
    \ip{\vect{\bhat}{i}}{\vect{\bhat}{j}} = \left\{
    \begin{array}{ll}
        0,                     & i\neq j\, , \\
        \|\vect{\bhat}{i}\|^2, & i=j\,.      \\
    \end{array}
    \right.
\end{equation}
For $0\leq t\leq n - m$ and $\vv v \in \RR^{n+1}$, the \emph{orthogonal projection}
$\awproj{t}(\vv v)$ is defined by
\[
    \awproj{t}: \RR^{n+1}\rightarrow \spann{b}{t}^\bot,\quad
    \vv v \mapsto \sum_{j=t}^{n - m}
    \frac{\ip{\vv v}{\vect{\bhat}{j}}}{\ip{\vect{\bhat}{j}}{\vect{\bhat}{j}}}
    \cdot \vect{\bhat}{j}
\]
and $\vect{\bhat}{t} = \awproj{t}(\vect{b}{t})$.

\bigskip

In each level $t$ of the backtracking algorithm,
$\vect{w}{t}  = \awproj{t}(\sum_{j=t}^{n-m} u_j\vect{b}{j})$
is the projection of the linear combination
of the already fixed variables $u_t$, $u_{t+1}$, $\ldots$, $u_{n-m}$
into the subspace of $\RR^{n+1}$ which is orthogonal to the linear span
$\langle b_0,\ldots,b_{t-1}\rangle$.

Starting from $\vect{w}{n-m+1}=\vv 0$,
$\vect{w}{t}$ can be computed iteratively from $\vect{w}{t+1}$ by
\[
    \vect{w}{t} = (\sum_{i=t}^{n-m}
    u_i \mu_{it})\vect{\bhat}{t} + \vect{w}{t+1}
\]
with Gram--Schmidt coefficients $\mu_{it}$.
In each level $t$, $n-m\geq t \geq 0$,
all possible integer values for the variable $u_t$ are tested.
The following two main tests allow to restrict the possible values of $u_t$.

\paragraph{First pruning condition.}
For all $j\leq t$, the vectors
$\vect{\bhat}{j}$ are orthogonal to $\vect{w}{t+1}$ and therefore
\[
    \|\vect{w}{t}\|_2^2 =
    (\sum_{i=t}^{n-m} u_i \mu_{it})^2\|\vect{\bhat}{t}\|_2^2 +
    \|\vect{w}{t+1}\|_2^2\;.
\]
That is, the sequence $\|\vect{w}{t}\|_2^2$ for $t=n-m, \ldots, 0$ is
increasing and allows a backtracking approach.
We notice that
$
    \vect{w}{0} = \sum_{j=0}^{n-m} u_j \vect{b}{j}
$
and
we can backtrack as soon as
\[
    \|\vect{w}{t}\|_2^2 > c := (n+1)\cdot \rmax^2\;.
\]
For fixed $u_{t+1}$, $\ldots$, $u_{n-m}$, this gives a bound for $u_t$:
\[
    (u_t + \sum_{i=t+1}^{n-m} u_i \mu_{it})^2
    \leq \frac{ c - \|\vect{w}{t+1}\|_2^2}{\|\vect{\bhat}{t}\|_2^2}\;.
\]

\paragraph{Second pruning condition.}
The second test is an adaption to the special situation
that we are searching for an integer linear combination
of the basis vectors which consists solely of components whose absolute value
is bounded by $\rmax$.
It is based on the following theorem by Ritter. 
\begin{theorem}[Ritter~\cite{Ritter:97}]\label{thm:ritter}
    If the given sequence of integers
    $u_t$, $u_{t+1}$, $\ldots$, $u_{n-m}\in \ZZ$ can be extended to
    $u_0$, $\ldots$, $u_t$, $\ldots$, $u_{n-m}\in \ZZ$ such that
    $\sum_{i=0}^{n-m} u_i\vect{b}{i}$ is a solution of \eqref{diophant},
    then for all $y_t$, $y_{t+1}$, $\ldots$, $y_{n-m}\in\RR$\/:
    \[
        |\sum_{i=t}^{n-m}
        y_i \|\vect{w}{i}\|_2^2| \leq \rmax \cdot
        \|\sum_{i=t}^{n-m} y_i\vect{w}{i} \|_1\;.
    \]
\end{theorem}
\noindent
We use this theorem in the enumeration algorithm in the following way.
Taking
$$(y_t,y_{t+1},\ldots,y_{n-m}) = (1,0,\ldots,0)$$
results in the test
\[
    \|\vect{w}{t}\|_2^2 \leq \rmax \|\vect{w}{t}\|_1\,. \label{aw:hoelder}
\]
If this inequality is violated for some vector $\vect{w}{t}=x\vect{\bhat}{t}+\vect{w}{t+1}$,
then it will also fail for all vectors of the form
$(x+r)\vect{\bhat}{t}+\vect{w}{t+1}$ with $r\in \ZZ$ and
$xr>0$.

\alglanguage{pseudocode}
\begin{algorithm}
    \caption{Lattice point enumeration}\label{algorithm1}
    Given the generator matrix~\eqref{lattice}
    of the lattice $\LL\subset\RR^{m+n+1}$ of rank $n+1$
    of problem \eqref{diophant},
    all nonzero vectors $\vv v\in \LL$ such that
    $\|\vv v\|_\infty\leq \rmax$ are determined.
    \begin{enumerate}
        \item
              Compute a reduced basis $\vect{b}{0},\vect{b}{1},\ldots,\vect{b}{n}$ of lattice $\LL$.
        \item
              Delete the unnecessary columns and rows of the generator matrix.
              The remaining basis $\vect{b}{0},\vect{b}{1},\ldots,\vect{b}{n-m}\subset\RR^{n+1}$ has rank $n-m+1$.
        \item
              Compute Gram--Schmidt orthogonalisation $\basisfull{\bhat}{n-m}$ together with
              Gram--Schmidt coefficients $\mu_{ij}$.
        \item Set $R := (n+1)\cdot\rmax^2$.
    \end{enumerate}

    \begin{algorithmic}[1]
        \Function{enum}{$t$, $\vv w'$}
        \State $\mbox{onedirection} \gets \mbox{false}$
        \State $y_t \gets \sum_{i=t+1}^{n-m} u_i \mu_{it}$
        \State $u_t \gets \lfloor-y_t\rceil$                 \label{test1} 
        \While{\mbox{true}}
        \State $\vv w \gets (\sum_{i=t}^{n-m} u_i \mu_{it})\vect{\bhat}{t} + \vv w'$
        \If{$\|\vv w\|_2^2 > R$}
        \Return {\ } \Comment{step back}
        \EndIf
        \If{$t>0$}
        \If{$\|\vv w\|_2^2 >  \rmax\cdot\|\vv w\|_1$} \label{test2} 
        \If{$\mbox{onedirection}$}
        \Return \Comment{step back}
        \Else
        \State next($u_t$)
        \Comment{select next $u_t$ in level $t$}
        \State $\mbox{onedirection} \gets \mbox{true}$
        \State \textbf{continue}
        \EndIf
        \Else
        \State enum($t-1$, $\vv w$)
        \Comment{step forward}
        \EndIf
        \Else
        \Comment{$t=0\; {}\rightarrow{}$ solution}
        \State output solution
        \EndIf
        \EndWhile
        \EndFunction
        \State enum($n-m, \vv 0$)
        \Comment{Start the recursive backtracking}
    \end{algorithmic}
\end{algorithm}

\begin{figure}[ht]
    \begin{center}
        \includegraphics{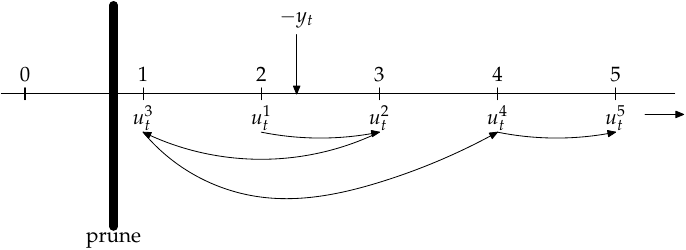}
    \end{center}
    \caption{Enumeration in level $t$ and pruning one direction after $u_t^3$}\label{fig:pruning}
\end{figure}

Algorithm~\ref{algorithm1} contains a high level description of the algorithm to solve \eqref{diophant},
where the procedure next() in Algorithm~\ref{algorithm1}
determines the next possible integer value of the variable $u_t$ in level $t$, compare
Figure~\ref{fig:pruning}.
Initially, when entering a new level $t$,
in line~\ref{test1} $u_t$ is set to be the closest integer value of $-y_t := -\sum_{i=t+1}^{n-m} u_i \mu_{it}$,
say $u_t^1$.
The next value $u_t^2$ of $u_t$ is the second closest integer to
$-y_t$ followed by $u_t^3$ and so forth. In other words, the values of
$u_t$ alternate with increasing distance around $-y_t$.

If the condition in line~\ref{test2} is true for the first time then we do one more regular call of the procedure next(),
i.e. $u_t$ is set to be the next closest integer to $-y_t$.
In the example of Fig.~\ref{fig:pruning} this happens when $u_t^4$ is determined.
After that, the enumeration proceeds only in this remaining direction, see
the computation of $u_t^5$ in Fig.~\ref{fig:pruning}.
Finally, when the condition  in line~\ref{test2} is true for the second time, the algorithm
steps back and increases the enumeration level.

\subsection{Recent improvements}

Since~\cite{wassermann:02} multiple improvements have been implemented
in the author's software \texttt{solvedi\-o\-phant}.

\paragraph{Speed improvements.}
Most loops to compute norms of vectors or involve vectors in general
have been sped up using SIMD instructions for x86 CPUs. At the time
of writing, AVX2 is supported. Alternatively, any BLAS library
like~\cite{openblas} might be used by setting a compiler flag.

\paragraph{Limited discrepancy search.}
\cite{10.1007/978-3-030-79987-8_2} discusses the use of
\emph{limited discrepancy search}~\cite{HarveyGinsberg:95} as an alternative
enumeration strategy to depth-first search.
This search paradigm proves to be valuable in some of the problems
where only one solution is needed like it is the case in (fMSP).


The enumeration algorithm of Algorithm \ref{algorithm1} performs depth-first search.
In particular, when entering enumeration level $t$, $u_t$ is chosen for
$\vv w := (\sum_{i=t}^{n-m} u_i \mu_{it})\vect{\bhat}{t} + \vv w'$ in line (7)
such that $\|\vv w\|_2$ is minimal among all choices for $u_t$.
In other words, the depth first search is organized using the heuristic
that choosing in each level the vector $\vv w$ such that $\|\vv w\|_2$ is minimal
will most probably lead to a solution vector.
However, it may be that this choice for $u_t$ in one of the first levels might
lead to no solution, but nevertheless the algorithm will enumerate a huge search tree below $u_t$.

This is a general problem of depth-first search.
In 1995, Harvey and Ginsberg~\cite{HarveyGinsberg:95} described an enumeration scheme
called \emph{limited discrepancy search} which aims to overcome this weakness.

Assume that a backtrack algorithm has to examine a search tree.
Each level corresponds to a variable and the algorithm has to assign a value to that variable,
followed by a test if this assignment might lead to a solution. If yes, we can proceed to the next level,
otherwise we have to assign a different value.
If we have tested all values, we have to step back to the previous level.
If values could be assigned to all variables, a solution has been found.

We assume that variable ordering is fixed and in each level of the backtracking
there exists a heuristic which determines the order in which the values are assigned
to the variable corresponding to that enumeration level.
A \emph{discrepancy} is defined as a deviation from the heuristic.

Harvey and Ginsberg suggest to enumerate the search tree in increasing number of
discrepancies.
In the first step, only the optimal choice in each level of the enumeration in Algorithm~\ref{algorithm1}
is assigned to the variables until there is a contradiction or a solution is found.
In the next step, all possible paths in the search tree with exactly one deviation (i.e. discrepancy)
from the heuristic are examined. After that, all paths in the search tree with two
deviations from the optimal choice are enumerated, and so forth.

The lattice point enumeration of Algorithm~\ref{algorithm1} and its value order heuristic show sometimes
dramatic improvements for finding the first solution in (fMSP) instances. It performs better than
depth-first search for all problems with $m\geq 7$, see Section~\ref{sec:results}.

\paragraph{Numerical stability.}
Since~\cite{MorelSV09,Schnorr10,Stehle10}, orthogonalization in the LLL algorithm is realized with Householder transformations
to improve numerical stability of the Gram--Schmidt orthogonalization.

Since~\cite{10.1007/978-3-030-79987-8_2}, the numerical stability of the lattice basis
reduction in Step 1 of the algorithm has been further enhanced as suggested in~\cite{OgitaRump:05}.
Even if numerical stability is not a problem for (fMSP) instances up to $m=15$, it is
important for combinatorial construction problems for which the number of variables exceeds $2000$.

\section{Results}\label{sec:results}

The algorithm described in the previous section has been implemented in
C using the AVX2 SIMD instruction set in the author's open source
software \texttt{solvediophant}.\footnote{Available from \href{https://github.com/alfredwassermann/solvediophant}{https://github.com/alfredwassermann/solvediophant}.}
We tested the QOBLIB market split instances with this program
on a computer with eight physical Intel Xeon E-2288G CPU (3.70GHz) processors and 64 GB memory, running Debian 12.
The tested version was \texttt{solvediophant} v1.0, using a single processor.

\begin{table}[ht]
    \caption{Computing results for instances from Quantum Optimization Benchmark Library~\cite{koch2025quantumoptimizationbenchmarklibrary}.}\label{table:timings}
    \begin{tabular}{lrrrr}
        \hline\noalign{\smallskip}
        Class          & \cite{KempkeKochZIB:25} (sec) & First (sec) & All (sec) & Number of solutions    \\
        \noalign{\smallskip}\hline\noalign{\smallskip}
        $(3,20,50)$    &                            & 0.00        & 0.00       & 1, 3, 1, 2             \\
        $(3,20,100)$   &                            & 0.00        & 0.00       & 1, 1, 1, 1             \\
        $(3,20,200)$   &                            & 0.00        & 0.00       & 1, 1, 1, 1             \\
        $(4,30,50)$    &                            & 0.00        & 0.00       & 1, 1, 2, 2             \\
        $(4,30,100)$   &                            & 0.00        & 0.00       & 1, 1, 2, 1             \\
        $(4,30,200)$   &                            & 0.00        & 0.00       & 1, 1, 1, 1             \\
        $(5,40,50)$    &                            & 0.00        & 0.00       & 23, 14, 16, 14         \\
        $(5,40,100)$   &                            & 0.00        & 0.00       & 2, 1, 2, 1             \\
        $(5,40,200)$   &                            & 0.00        & 0.00       & 1, 1, 1, 1             \\
        $(6,50,50)$    &                            & 0.00        & 2.00       & 45, 37, 53, 40         \\
        $(6,50,100)$   &                            & 0.75        & 1.50       & 1, 1, 1, 1             \\
        $(6,50,200)$   &                            & 1.00        & 1.75       & 1, 1, 1, 1             \\
        $(7,60,50)$    & 0.37                       & 0.75        & 9.25       & 25, 180, 306, 178      \\
        $(7,60,100)$   & 1.20                       & 1.50        & 5.50       & 1, 4, 2, 1             \\
        $(7,60,200)$   & 2.02                       & 3.00        & 4.75       & 1, 1, 1, 1             \\
        $(8,70,50)$    & 1.25                       & 3.75        & 505.50     & 1265, 1066, 752, 943   \\
        $(8,70,100)$   & 7.84                       & 6.25        & 93.00      & 3, 5, 3, 4             \\
        $(8,70,200)$   & 17.80                      & 6.00        & 24.50      & 1, 1, 1, 1             \\
        $(9,80,50)$    & 16.03                         & 1.25        & 36,825.00 & 4497, 3720, 3135, 3247 \\
        $(9,80,100)$   & 236.07                        & 10.00       & 4,203.75  & 7, 7, 6, 7             \\
        $(9,80,200)$   & 520.76                     & 48.50       & 824.75     & 1, 1, 1, 1             \\
        $(10,90,50)$   & 167.83                     & 6.75        & --         & --                     \\
        $(10,90,100)$  & 66,636.22                     & 1,169.25    & --        & --                     \\
        $(10,90,200)$  &                               & 13,881.50   & 61,186.50 & 1, 1, 1, 1             \\
        $(11,100,50)$  & 41,399.39                     & 24.00       & --        & --                     \\
        $(11,100,100)$ &                               & 70,421.50   & --        & --                     \\
        $(11,100,200)$ &                               & 453,668.50  & --        & --                     \\
        $(12,110,50)$  &                            & 617.00      & --         & --                     \\
        $(13,120,50)$  &                               & 5,365.00    & --        & --                     \\
        $(14,130,50)$  &                               & 140,823.00  & --        & --                     \\  
        $(15,140,50)$  &                            & --          & --         & --                     \\
        \noalign{\smallskip}\hline
    \end{tabular}
\end{table}

As described in~\cite[Section 4.1]{koch2025quantumoptimizationbenchmarklibrary},
there are four instances for each $m\in\{3,\ldots,15\}$ and $D\in\{50, 100, 200\}$
and all instances are feasible by construction.

\cite{KempkeKochZIB:25} reports benchmark results for these problem sets with their implementation of
the Schroe\-ppel--Shamir algorithm
as well as with the integer linear programming formulation of the problem using Gurobi.
Since the full number of solutions for an instance is never mentioned in that paper
it seems that also in that paper an instance of (fMSP) is regarded as \emph{being solved} as soon as one solution has been found.
Here, we report the run times to solve the instances as well as the run times to find all solutions
of (fMSP) instances.

All instances for $m\in\{3,\ldots,15\}$ and $D\in\{50, 100, 200\}$ have been tested,
Table~\ref{table:timings} shows the results.
Column ``Class'' contains the problem class, i.e.
$(m, 10(m-1), D)$. For each class, all four instances have been tested and the average of the run times is listed.

Column ``\cite{KempkeKochZIB:25}'' lists the average running time in seconds to find one solution with the Schroep\-pel--Shamir algorithm as
reported by~\cite{KempkeKochZIB:25} on a computer using a GPU. The authors published results for instances with
$m\in\{7,\ldots,11\}$. They give also results for Gurobi for these instances, but it was always
slower than Schroeppel--Shamir or could not solve the instances at all.
In the table, an empty entry means that in~\cite{KempkeKochZIB:25} no data has been reported for this class.

Column ``First'' contains the average time in seconds 
that \texttt{solvediophant} v1.0
needed for the four instances until the first solution has been found, i.e. the instance has been solved.
If an entry is $0.00$, the average run time was less than 0.5 seconds.
An entry ``--'' means that the enumeration could not be completed in several days.

Column ``All'' lists the average run time in seconds for the four instances
that the algorithm needed to enumerate all solutions.
If an entry is $0.00$, the average run time was less than 0.5 seconds.
An entry ``--'' means that the enumeration could not be completed in several days.

Any run time given in columns ``First'' and ``All'' comprises the time for lattice basis reduction
and enumeration. No preprocessing of the instances had been performed.
The correctness of all solutions was tested with a completely independent Python program. The run time for these
tests is negligible.

The number of feasible solutions of each instance
is given in column ``Number of solutions''. Here, the ordering is according
to the lexicographic ordering of the file names. The exact number of solutions of the various instances
had not been reported previously and seems to be new.

Table~\ref{table:timings} shows that there exist many solutions of (fMSP) if $D=50$ and $m\geq 5$.
On the other hand, all instances for $D=200$ and $m\leq 11$ have exactly one solution. If one
assumes that the solutions are distributed uniformly in the leaves of the search tree, 
it is no surprise that the search for the first solution is significantly faster
for smaller values of $D$ than for larger values.

The data for column ``First'' was determined by
starting two programs simultaneously on two cores for each instance.
One program used enumeration by limited discrepancy search (LDS) and the other program
used enumeration by depth first search (DFS). For each instance, the time of the program that terminated first was taken into account.
It turned out that for all instances with $m\geq7$, LDS was the faster enumeration strategy, beating DFS.
For column ``All'', only DFS enumeration was used, since its overall bookkeeping overhead is less than for LDS.

\section{Discussion and future work}

This note shows that algorithms based on lattice basis reduction and lattice enumeration
seem to be good candidates to solve the market split problem.

If only one solution is needed, it is promising to  
use \emph{limited discrepancy search} in the second step of the algorithm.
This is the case for QOBLIB since all market split instances are feasible.
For all classes with non-negligible run time
limited discrepancy search was faster than depth-first search in finding the first solution.
This might be different for market split feasibility instances in the original formulation from~\cite{CornuejolsDawande:98}:
the administrative overhead for limited discrepancy search
is always larger than that of depth-first search (beside border cases)
and if an instance is infeasible, the algorithm has to exhaustively run through the whole
search tree.

The proposed algorithm can be parallelized by partitioning the work and distributing the parts of
the search tree to several computer nodes, see e.g.~\cite{bettenwassermann96}.
A different promising strategy for parallelization of the problem
is to start the algorithm simultaneously
on many computer nodes with the
order of the input lattice basis \eqref{lattice} shuffled randomly, see~\cite{aono:2025}
for the theoretical background.


\section*{Statements and Declarations}

\medskip\noindent
The datasets generated during and analysed during the current study are from the QOBLIB repository,
available at
\\
\centerline{\href{https://github.com/ZIB-AOPT/QOBLIB}{https://github.com/ZIB-AOPT/QOBLIB}}
\\
and
\\
\centerline{\href{https://git.zib.de/qopt/qoblib-quantum-optimization-benchmarking-library}{https://git.zib.de/qopt/qoblib-quantum-optimization-benchmarking-library}.}

\medskip\noindent
\texttt{solvediophant} is available from 
\\
\centerline{\href{https://github.com/alfredwassermann/solvediophant}{https://github.com/alfredwassermann/solvediophant}.}



\end{document}